\documentclass[12pt,reqno,twoside]{amsart}

\usepackage{epsfig}
\usepackage{color}
\usepackage[english]{babel}
\usepackage{amsmath,amssymb,enumerate,amsthm}
\usepackage[latin1]{inputenc}
\usepackage[T1]{fontenc}
\usepackage{cite}
\usepackage{latexsym}
\usepackage{setspace}
\usepackage{bm} 
\usepackage{comment}
\usepackage{array,graphicx,caption}
\newtheorem{theorem}{Theorem}
\newtheorem*{theorem*}{Theorem}

\newtheorem*{conjecture*}{Conjecture}

\newtheorem*{definition*}{Definition} 
\newtheorem{lemma}{Lemma}
\newtheorem*{lemma*}{Lemma}

\newtheorem*{corollary*}{Corollary}

\newtheorem*{proposition*}{Proposition}

\begin{document}

\begin{Large}
\centerline{\bf
On badly approximable subspaces}
\vskip+0.5cm
\centerline{Natalia Dyakova\footnote{Research is supported by RFBR grant no. 18-01-00886a.}}
\end{Large}

\vskip+1.5cm
\section{Introduction}
In a remarkable paper \cite{Kle1}  devoted  to the 
dynamical approach to Diophantine approximations on submanifolds
 D. Kleinbock 
among other results observed the following phenomenon. If an affine subspace 
of Euclidean space contains a vector which is "not very well approximable" by rational vectors, then 
almost all
  vectors in this affine subspace
  are "not very well approximable". A similar result holds for manifolds 
  and was generalized in different directions (see \cite{Kle2}, \cite{Kle3} and
   \cite{Zha}). 
In \cite{Moshch} N. Moshchevitin gave a simple proof of the result by Kleinbock.
In the present paper we prove a general result about badly appriximable susbpaces
of a certain linear subspace of $\mathbb{R}^d$.  The result by Moshchevitin from \cite{Moshch} is a particular case of our main theorem.

\subsection{Badly approximable  matrices and linear subspaces}

Let
\begin{gather}\label{Theta}
\Theta =
\left(
\begin{array}{ccc}
\theta_{1,1}&...& \theta_{1,m}\cr
\vdots&\vdots&\vdots \cr
\theta_{n,1}&...& \theta_{n,m}\cr
\end{array}
\right)
\end{gather}
be a real matrix and
$\psi(t)$ be a real valued function decreasing to zero as $ t \to \infty$.
Matrix $\Theta$ is called $\psi$-badly approximable if for any integer vector $ \pmb{x}= (x_1,...,x_m) \in \mathbb{Z}^m$
one has
$$
\max_{1\le j\le n} ||\theta_{j,1}x_1+...+\theta_{j,m}x_m||
\ge \psi (|\pmb{x}|).
$$
Here $||\xi || = \min_{x\in \mathbb{Z}}  |\xi - x|$
is the distance from real $\xi$ to the nearest integer and
$|\pmb{\xi}| = \max_{1\le i \le m} |\xi_i|$
stands for the sup-norm of the vector $\pmb{\xi} =(\xi_1,...,\xi_m)\in \mathbb{R}^m$.
The irrationality measure function $\psi_\Theta (t)$ associated with matrix $\Theta$ is defined as
\begin{equation}\label{psiTheta}
\psi_\Theta (t) =
\min_{ \substack{\pmb{x} \in \mathbb{Z}^m\\
0<
        |\pmb{x}|\le t}}
\,\,\,
\max_{1\le j\le n} ||\theta_{j,1}x_1+...+\theta_{j,m}x_m||,\,\,\,\,\, t \ge 1.
\end{equation}
Then matrix $\Theta$ is $\psi$-badly approximable if and only if
$$
\psi_\Theta (t) \ge \psi (t),\,\,\,\,\, \forall\, t \ge 1.
$$
We should note that for any positive $\varepsilon$  almost all matrices $\Theta$ are $\psi$-badly approximable 
(with respect to Lebesgue  $(m\times n)$-dimensional measure)
with
$\psi (t) = \rho t^{-m/n -\varepsilon}$, for some $\rho>0$. Moreover the set
$$
\{\Theta:\,\, \Theta \,\,\text{is }\, \rho t^{-m/n}\text{-badly approximable for some }\,\, \rho >0\}
$$
is of zero measure  and has full Hausdorff dimension
(see  \cite{Schm,Schm1}).

However, it will be more convenient for us to consider a little bit different irrationality measure function.

Let $d=m+n$ and $ \mathbb{R}^{d}$ be $d$-dimensional Euclidean space with coordinates 
$$\pmb{z} = (\pmb{x},\pmb{y})=(x_1,...,x_m,y_1,...,y_n) .$$

For a proper linear subspace $\mathfrak{B} \subset \mathbb{R}^d$  we consider the function

\begin{equation}\label{psiB}
\psi_{\mathfrak{B}}(t) = \min\limits_{\substack{\bm{z}\in \mathbb{Z}^d \\ 1 \le |\bm{z}|<t}}{\rm dist} (\bm{z},\mathfrak{B}),
\end{equation}
where ${\rm dist} (A,B)$ denotes the Euclidean distance between the sets $A$ and $B$.

 We define a proper linear subspace $\mathfrak{B} \subset \mathbb{R}^d$ to be $\psi$-badly approximable if
 $$
 \psi_{\mathfrak{B}}(t) \ge\psi(t),\quad \forall t \ge 1.$$

Now, for the matrix $\Theta$ we consider the linear subspace
\begin{equation}\label{LTheta}
\mathfrak{L}_\Theta =
\{
\pmb{z} =(x_1,...,x_m,y_1,...,y_n) \in
\mathbb{R}^{m+n}:\,\,
\pmb{y} = \Theta \pmb{x}\}.
\end{equation}
Here we should note that
the projection $z=(x_1,...,x_m,y_1,...,y_n)\mapsto  (x_1,\ldots, x_m)$
gives a bijection between $\mathfrak{L}_{\Theta}$ and coordinate subspace $\{y_1 = \ldots = y_n = 0\}$.
It is clear that
$$
\kappa_1 \psi_\Theta (\kappa_2 t) \le \psi_{\mathfrak{L}_{\Theta}} (t) \le \kappa_3 \psi_\Theta (\kappa_4 t)
$$
with some positive $\kappa_j$ depending on $\Theta.$

So, if the subspace $\mathfrak{B} = \mathfrak{L}_\Theta$ is of the form \eqref{LTheta}, then the property of $\mathfrak{B}$ to be $\psi-$badly approximable means that $$
\psi_\Theta \ge \kappa' \psi(\kappa'' t)
$$
 with certain  positive $\kappa', \kappa''$. Similarly $\Theta$ is $\psi$-badly approximable matrix means that
 $$
\psi_{\mathfrak{L}_\Theta} \ge \kappa^* \psi(\kappa^{**} t)
$$
 with certain  positive $\kappa^*, \kappa^{**}$.

\subsection{Diophantine exponent}
For a real matrix $\Theta$ we define the {\it Diophantine exponent} as 
$$
\omega(\Theta) = \sup \{ \tau : \liminf_{t \rightarrow \infty} t^{\tau} \psi_{\Theta} (t) < + \infty
\}.
$$
From the discussion of previous subsection it is clear that 
 if $\mathfrak{B} = \mathfrak{L}_\Theta$ then 
 $$\omega(\mathfrak{B}) := \sup \{ \gamma: \liminf_{t \rightarrow \infty} t^\gamma \cdot \psi_{\mathfrak{B}}(t) < \infty\} = \omega(\Theta).$$ 

If $\mathfrak{B} = \mathfrak{L}_\Theta$, so the equality $\omega = \omega(\mathfrak{B})$ 
means that 
 $\forall \tau < \omega(\mathfrak{B})$ there exists infinitely many $\pmb{z} = (x_1, \ldots, x_m,y_1, \ldots, y_n) \in \mathbb{Z}^d$ such that
 $$
\max_{1\le j\le n} |\theta_{j,1}x_1+...\theta_{j,m}x_m - y_j|
< (\max_{1\le i\le m} |x_i|)^{-\tau},
 $$
 or $$
 \omega(\mathfrak{B}) = \inf \{\gamma : \Theta \,\, \text{is $\rho t^{-\gamma}$-badly approximable for some $\rho >0$}\,\,
 \}.
 $$

\subsection{Parametrization of subspaces} 

Our results deal with the following situation. Given a linear $a$-dimensional subspace $\mathfrak{A} \subset \mathbb{R}^d$ we want to get a statement about {\it almost all} $c$-dimensional linear subspaces $\mathfrak{C} \subset \mathfrak{A}$. So for a fixed $\mathfrak{A}$ we should give a parametrization of all linear  $c$-dimensional subspaces of $\mathfrak{A}$.

First of all we choose Euclidean coordinates 
$$
\pmb{w} = (\pmb{u},\pmb{v}), \,\,\,
\pmb{u} = (u_1,...,u_a), \,\,\, \pmb{v} = (v_1,...,v_{d-a})
$$
in $\mathbb{R}^d$ in such a way that
the subspace $\mathfrak{A} $ is defined by
$$
v_1=...=v_{d-a} = 0.
$$
We shall use these coordinates everywhere in our proofs.
In these coordinates the  lattice of integer points in $\mathbb{R}^d$ will be denoted as $\Lambda$. It is clear that its fundamental volume is equal to one.

Given a $c$-dimensional linear subspace $\mathfrak{C} \subset \mathfrak{A}$. There exists a coordinate subspace 
$$\mathbb{R}^c =\{u_{i_1}= \ldots= u_{i_{a-c}} =0\}$$ 
of $\mathfrak{A}$ such that the projection $\mathfrak{C} \rightarrow \mathbb{R}^c$ is bijective. So we restrict ourselves with the case when 
\begin{equation}\label{ce}
\mathfrak{C} =  \mathfrak{C}  (\Theta) = \{\bm{u} =( \bm{x}, \bm{y}  ):\,\, \bm{x}\in \mathbb{R}^{c}, \bm{y} \in \mathbb{R}^{a-c} \,\,\text{such that}\,\, \bm{y} = \Theta \bm{x} \}
\end{equation}
 where $\Theta = (\theta_{ij})  $  is a $(a-c)\times c$-matrix with real entries, as it was done for example in \cite{GG}.
 If $\mathfrak{C} $ satisfies (\ref{ce}) for the corresponding matrix $\Theta$ we write
 $\Theta = \Theta (\mathfrak{C} )$.


So in the present paper, when we speak about almost all $c$-dimensional linear subspace $\mathfrak{C} \subset \mathfrak{A}$ we mean almost all  $(a-c)\times c$-matrices with respect to Lebesgue measure in $\mathbb{R}^{(a-c) \times c}$. An alternative approach to consider the variety of all $c$-dimensional subspaces of $\mathfrak{A}$
is connected with Grassman coordinates (see \cite{Sccc}, Ch. 1, \S 5). 

\subsection{The main result for Diophantine exponents}
Here we formulate our main  results in terms of Diophantine exponents of linear subspaces of $\mathbb{R}^d$.
 \begin{theorem}\label{th2}
 Let $\mathfrak{A}$ be $a$-dimensional linear subspace of $\mathbb{R}^d$. Suppose that $\mathfrak{B}$ is $b$-dimensional linear subspace of $\mathfrak{A}$ with $\omega(\mathfrak{B}) = \omega$.

1)  If $c < b$, then $$\omega(\mathfrak{C}) \le \omega +\dfrac{(\omega +1) \cdot (c-b)}{d-c} =  \dfrac{\omega \cdot (d-b) +c-b}{d-c}$$ for almost all $c$-dimensional linear subspaces $\mathfrak{C} \subset \mathfrak{A}$.

2) If $c \ge b$, then $$\omega(\mathfrak{C}) \le \omega +\dfrac{(\omega +1) \cdot (c-b)}{a-c} =  \dfrac{\omega  \cdot (a-b) +c-b}{a-c}$$ for almost all $c$-dimensional linear subspaces $\mathfrak{C} \subset \mathfrak{A}$.
\end{theorem}
Theorem \ref{th2} means that if there exists $b$-dimensional linear subspace $\mathfrak{B} \subset \mathfrak{A}$ which doesn't have good rational approximations then almost all $c$-dimensional subspaces $\mathfrak{C} \subset \mathfrak{A}$  don't have good rational approximations  either.
 Theorem \ref{th2}  is a particular case of more general result (Theorem \ref{th1}), which will be formulated in the next section.
\subsection{An example}

We do not  want to give here examples of application our results for $c=1$ as this case was completely considered in \cite{Moshch}. 
To illustrate our theorem we consider the following example.

Let us deal with $2$-dimensional badly approximable subspace in $\mathbb{R}^4$. It is clear that there exists $T^{-\beta}$-badly approximable $2$-dimensional subspace with $\beta = 1.$

However if we would like to find a badly approximable two-dimensional subspace in a three-dimensional subspace $\mathfrak{A} \subset \mathbb{R}^4$, dim $\mathfrak{A} = 3$, it may happen that we be able to find $T^{-\beta}$ badly approximable $2$-dimensional subspace with $\beta = 2$ only.
Indeed, this happens if we take $\mathfrak{A} \subset \mathbb{R}^4$ to be a completely rational $3$-dimensional linear subspace.

Our  Theorem 1   gives the following results. Suppose that $\mathfrak{A} \subset \mathbb{R}^4$ is a $3$-dimensional linear subspace of $\mathbb{R}^4$ such that there exists a vector $\pmb{w} = (1, w_1, w_2, w_3) \in \mathfrak{A}$ satisfying
 $$\max_{1\le i\le 3} || w_i q || \ge \dfrac{c}{q^{1/3}} $$
 for all $ q\in \mathbb{Z}_+$ with some positive $c$
(so, $1$-dimensional linear subspace spanned by $\pmb{w}$ is $cT^{-\frac{1}{3}}$-badly approximable).
Then almost all $2$-dimensional linear subspace $\mathfrak{C} \subset \mathfrak{A}$ are $\rho T^{-\beta}$-badly approximable with any $\beta >\frac{5}{3}$ and some positive $\rho$. 
\section{General result}
\subsection{Functions and examples}

Let $\psi (T)$ and $\varphi(T)$, $T \ge 1 $ be two positive valued functions monotonically decreasing to zero as $T \rightarrow \infty$
{
satisfying special technical conditions.
}
Here we deal with $\psi$-badly approximable $b$-dimensional linear subspace $\mathfrak{B} \subset \mathbb{R}^d$. We may suppose that
$$\psi(T) \le T^{-\frac{b}{d-b}}.$$ 
{
For these functions we define quantities  $\mu_T = 0 $ for $ 0\le T <1 $ and
$$
  \mu_T = \left( \frac{T}{\psi(T)}\right)^{a-b} \cdot\max \left\{1, \left(\frac{\varphi(T) }{\psi (T)}\right)^{d-a}\right\}
  =
  $$
   \begin{equation}\label{mu_T}
   =
  \begin{cases}
   \left( \frac{T}{\psi(T)}\right)^{a-b},\,\,\,\,\varphi(T) \le \psi (T)\cr
   \frac{T^{a-b}\cdot\varphi^{d-a}(T)}{\psi^{d-b}(T)},\,\,\,\,
   \varphi(T) \ge \psi (T)
   \end{cases}
    \end{equation} 
    for $ T\ge 1$.
    For integer values of $T$
    we consider the quantities
    $$
    M_0 = 0,\,\,\,\text{and}\,\,\,
     M_T = \sum_{j=0}^{[\log_2T]} \mu_{T/2^j},\,\,\,\,\,
  \lambda_T = M_{T} - M_{T-1}, \,\,\, T \ge 1,
  $$
  so
  \begin{equation}\label{maa}
  M_T= \sum_{j=1}^T \lambda_j.
  \end{equation}
 }
\begin{theorem}\label{th1}
{
Let  all the functions satisfy conditions of this section and  the series 
\begin{equation}\label{series}
 \sum_{T=1}^{+\infty} \lambda_T  \cdot \dfrac{\varphi^{a-c}(T)}{T^{a-c}}
\end{equation}}
converge.

Suppose that $\mathfrak{B}$ is $b$-dimensional $\psi$-badly approximable linear subspace. 
Then for each $a$-dimensional linear subspace $\mathfrak{A} \supset \mathfrak{B}$
almost all  $c$-dimensional linear subspaces $\mathfrak{C} \subset \mathfrak{A}$ are $\varphi$-badly approximable subspaces.
\end{theorem}

{\bf Remark.} \,
For $ c = 1$ the result of Theorem 2 {is similar to}  the main result form the paper \cite{Moshch}.


To illustrate our Theorem \ref{th1} we consider the function 
$
\psi (T) = T^{-\beta} \log^B T.
$
Then we can choose 
$
\varphi  (T) = T^{-\gamma} \log^C T
$
with

 \vskip+0.3cm
$ \left[
 \begin{gathered}
 \gamma > \dfrac{\beta \cdot (d-b) +c-b}{d-c} \quad {\text where} \ C \ {\text and} \ B \ {\text are} \ {\text arbitrary}\\
 or\\
 \gamma = \dfrac{\beta \cdot (d-b) +c-b}{d-c} \quad {\text and} \ C (d-c) - B (d-b) < -1,\\
 \end{gathered}
 \right.$
 \vskip+0.3cm
 and $ \gamma <\beta$ for the case when $c < b$, or
 
  \vskip+0.3cm
 $ \left[
 \begin{gathered}
\gamma >\dfrac{\beta \cdot (a-b) +c-b}{a-c}\quad {\text where} \ C \ {\text and} \ B \ {\text are} \ {\text arbitrary}\\
 \\
or\\
 \\
\gamma = \dfrac{\beta \cdot (a-b) +c-b}{a-c} \quad {\text and} \ C (a-c) - B (a-b) < -1,\\
 \end{gathered}
 \right.$
  \vskip+0.3cm
 
 for the case when $c \ge b$
 (here we should note that under these inequalities we always have $ \gamma \ge \beta$).

  In both cases the series (\ref{series}) converges.
  Indeed, consider the case $c<b$. Then  we have $ \varphi (T) >\psi (T)$ and
 $$
 \mu_T = T^f\log^g T\,\,\,
 \text{with}\,\,\,
 \begin{cases}
 f= a-b-\gamma(d-a)+\beta (d-b)\cr
 g = C(d-a)-B(d-b)
 \end{cases},
 $$
 $$
 M_T = 
 \sum_{j=0}^{[\log_2T]}
 \left(\frac{T}{2^j}\right)^f\log^g \left(\frac{T}{2^j}\right),
 $$
  $$
  \lambda_T =
  M_{T}- M_{T-1} =
  $$
  $$
  =
   \sum_{j=0}^{[\log_2(T-1)]}\left(
 \left(\frac{T}{2^j}\right)^f\log^g \left(\frac{T}{2^j}\right)
 -
 \left(\frac{T-1}{2^j}\right)^f\log^g \left(\frac{T-1}{2^j}\right)
 \right)+ O(1) =
 $$
  $$
  =
 O(T^{f-1}\log^gT).
 $$
 So  the summands from the series  (\ref{series}) are of the form
 $$
 O(
 T^{f_1}\log^{g_1}T)\,\,\,
 \text{with}\,\,\,
  \begin{cases}
 f_1= 
 -\gamma (d-c)+\beta(d-b)+c-b-1\le -1\cr
 g_1 = C(d-c)-B(d-b)<-1
 \end{cases}.
 $$
 The case $c \ge b$ is similar.

We  see that for the functions $\varphi (T) $ and $\psi (T)$  under consideration the conclusion of the Theorem \ref{th1} is valid. In particular the choice of 
$$\psi(T) = t^{-\omega(\mathfrak{B}) - \varepsilon}, \quad \varepsilon >0$$
shows that the Theorem \ref{th2} from the previous section follows from Theorem~\ref{th1}. 

\section{Proofs}

\subsection{Simple lemmas}

Consider the sets
$$
\Omega_T = \{ \bm{w} = (w_1, \ldots,  w_d) \in \mathbb{R}^d:|\bm{w}|\le T , \inf\limits_{\bm{b} \in \mathfrak{B}}  
{\rm dist} \,(
\bm{w} , \bm{b})\le  \psi(T)  \}
$$
and
$$
\Pi_T = \{\bm{w} = (w_1, \ldots,  w_d) \in \mathbb{R}^d : |\bm{w}| \le T, \inf\limits_{\bm{a} \in \mathfrak{A}} {\rm dist} \,(
\bm{w} , \bm{a})\le \varphi(T)\}.
$$
\begin{lemma}\label{1}
Let $N_T$ be the minimal number of points $\bm{c}_j, 1\le j \le N_T$ such that
$$
\Pi_T \subset \bigcup_{j=1}^{N_T} \left(\frac{1}{2} \Omega_T +\bm{c}_j\right).
$$
Then {
\begin{equation}\label{ente}
N_T \ll \mu_T,
\end{equation}}
where $\mu_T$ is defined in \eqref{mu_T}.
\end{lemma}

{
{\bf Remark.
}
The same bound (\ref{ente}) holds for the number  of copies of the set $\Omega_T $ which is sufficient to cover any shifted set
$ \Pi_T +\pmb{z}$. So for the number  $N_T'$ of copies of the set $\Omega_T $ which is sufficient to cover $\Pi_{2T}$  the bound (\ref{ente}) holds also.}

\textbf{Proof of Lemma 1.}
As $\Pi_T, \Omega_T \subset E_T = \{\pmb{x}: |\pmb{x}|< T\}$, by $O\left(\left(\dfrac{T}{\psi (T)}\right)^{a-b}\right)$ shifted sets $\frac{1}{2} \Omega_T +\bm{c}_j$ we are able to cover $\dfrac{\psi(T)}{4}$ neighbourhood of $\mathfrak{A}\cap E_T.$

So, in the case $\varphi (T) \le \psi (T) $ by $O\left(\left(\dfrac{T}{\psi (T)}\right)^{a-b}\right)$ shifted sets $\frac{1}{2} \Omega_T +\bm{c}_j$ we will cover the whole $\Pi_T$. So, $$
N_T \ll  \left(\frac{T}{\psi(T)}\right)^{a-b}. 
$$ 
and in this case lemma is proved.
In the case $\varphi (T) > \psi (T)$, to cover $\varphi(T)$-neighbourhood of $\mathfrak{A}\cap E_T$ we need to cover $O\left(\left(\dfrac{\varphi(T)}{\psi(T)}\right)^{d-a}\right)$ shifted copies of $\psi(T)$-neighbourhood of $\mathfrak{A}\cap E_T$. So in this case
$$
N_T \ll  \frac{T^{a-b}\cdot\varphi^{d-a}(T)}{\psi^{d-b}(T)}. 
$$
and everything  is proved.
$
\square
$

\begin{lemma}\label{2}
Suppose that $\mathfrak{B}$ is $\psi$-badly approximable. Then for any $T>1$ and for any $\bm{c} \in \mathbb{R}^d$ 
 the set $\frac{1}{2} \times \Omega_T + \bm{c}$ consists of not more than one point from the lattice $\Lambda$ defined in Subsection 1.3.

\end{lemma}

\textbf{Proof.}
Since $\mathfrak{B}$ is a $\psi-$badly approximable subspace it is clear that for any $T$
$$
\Omega_T \cap \Lambda = \{ 0\}.
$$
Consider any translation of the  $1/2-$contracted set
$$
\frac{1}{2} \times \Omega_T + \bm{c},\qquad \bm{c} \in \mathbb{R}^d. 
$$
If two different integer points $\bm{w}, \bm{t}$ belong to the same set of this form, then $$
\bm{0} \neq \bm{w}-\bm{t} \in \Omega_T.
$$
But this is not possible.
$\square$

\subsection{Proof of Theorem \ref{th1}}

As was explained in the section 1.4 we are dealing with $a$-dimensional space $\mathfrak{A}$ with coordinates
$$
x_1, \ldots, x_c, y_1, \ldots, y_{a-c}.
$$
It is clear that any $c$-dimensional subspace under consideration $\mathfrak{C} \subset \mathfrak{A}$ can be defined by the system of $a-c$ equations

$$
\left\{
\begin{array}{rcl}
 \theta_{1,1} x_1+\ldots+ \theta_{1,c} x_c = y_1 \\
 \ldots \quad \ldots \quad \ldots \\
		\theta_{a-c, 1} x_1 + \theta_{a-c, c} x_c = y_{a-c}.\\
		\end{array}
\right.
$$
We take $R >1$ and suppose that 
\begin{equation}\label{tt}
\max_{i,j}| \theta_{i,j} | \le R.
\end{equation}
Consider $(a-1)$-dimensional subspaces $$
\mathfrak{C}_j =\{  (\pmb{x,y}) \in \mathfrak{A}  : \theta_{j,1} x_1+\ldots+\theta_{j,c} x_c = y_j\}.
$$
Then $$
\mathfrak{C} = \bigcap_{1\le j \le a-c}  \mathfrak{C}_j.$$

Let $\bm{z} \in \mathbb{Z}^d$ be an integer point.
 Denote by $\bm{z}^* = \bm{z}|_{\mathfrak{A}}$ the orthogonal projection of $\bm{z}$ onto $\mathfrak{A}$.
Consider the neighbourhood $$ U_{\varphi (|\bm{z}|)} (\bm{z}) = \{ \pmb{z'} \in \mathfrak{A} : \mbox{dist} (\pmb{z'},\pmb{z}) < \varphi (|\bm{z}|)\}.$$
We are interested in those $\pmb{z}$ for which  $$\mathfrak{A} \cap U_{\varphi ( |\bm{z}|)} (\bm{z}) \neq \emptyset.$$ 
These $\bm{z}$ are close $\mathfrak{A}$ and for them $|\bm{z}| \asymp |\bm{z}^*| $.

Consider the set of $c$-dimensional  subspaces 
$$
L_{\bm{z}^*} = \{
\mathfrak{C} =
\mathfrak{C} (\Theta)
\subset \mathfrak{A} :   \,\,\,   \Theta \,\,\text{satisfies (\ref{tt})},   \,\,\, \mbox{dist} (\mathfrak{C}, \bm{z}^*)
<  \varepsilon \}.$$
Our parametrization (see Section 1.4) allows to consider 
$$\Theta (L_{\bm{z}^*}) = \{ \Theta:\,\, \Theta = \Theta ( \mathfrak{C}), \,\, \mathfrak{C}\in  L_{\bm{z}^*} \}
$$ 
as a subsets of $\mathbb{R}^{(a-c)\times c}$.
The condition $\mbox{dist} (\mathfrak{C}, \bm{z}^*)
< \varepsilon$ leads to
$$
\begin{array}{rcl}
\mbox{dist} (\mathfrak{C}_i, \bm{z}^*) = \dfrac{| \theta_{i,1} z_1^{*}+\ldots+\theta_{i,c} z_c^{*} - z_{c+i}^{*} |}{\sqrt{\theta^2_{i,1}+\ldots+ \theta^2_{i,c}+1}} <  \varepsilon\\
		\end{array}
\qquad \forall  i = 1, \ldots, a-c.
$$

So 
$$
\Theta (L_{\bm{z}^*}) \subset \widehat{L}_{\bm{z}^*}=
\widehat{L}_{\bm{z}^*} (R, \varepsilon)=
$$
$$= \{ \Theta \in \mathbb{R}^{(a-c)\times c} : 
\mbox{dist} (\mathfrak{C}_i, \bm{z}^*) = | \sum_{j=1}^c \theta_{ij} z_j^{*} - z_{c+i}^{*} | < \sqrt{c+1} R \varepsilon 
$$
$$
\forall  i = 1, \ldots, a-c\\
		\}.
		$$

 Now we consider $\theta_{ij}$ as a variables. Then equations 
 $$
  \theta_{i,1} z_1^{*}+\ldots+\theta_{i,c} z_c^{*} - z_{c+i}^{*} =0\\
\qquad \forall  i = 1, \ldots, a-c.
$$
define an affine $(a-c)\cdot (c-1)$-dimensional subspace
$$\{ \Theta \in \mathbb{R}^{(a-c)\times c} : 
\begin{array}{rcl}
\mbox{dist} (\mathfrak{C}_i, \bm{z}^*) = 0 \quad\forall  i = 1, \ldots, a-c\\
		\end{array} \}.$$

 So the set  $L_{\bm{z}^*}$ is a neighbourhood of this affine subspace.
 
Now for $(a-c)\times c$-dimensional Lebesgue measure we have
\begin{equation}\label{varphi}
\mu(\widehat{L}_{\bm{z}^*}) \ll \left( \dfrac{R\varepsilon}{\max_{1\le i\le c}|z_i^*|}\right)^{a-c} R^c \ll
\dfrac{R^{a+2}\varepsilon^{a-c}}{|\pmb{z}|^{a-c}}
.
\end{equation}
In the last inequality we use  the upper bound 
$|z_{c+i}^*| \ll R (|z_1^*|+...+|z_c^*|)$  for all $i\in \{1,...,a-c\}$ which follows form
(\ref{tt}).

 

Now we suppose that $\mathfrak{C} \subset \mathfrak{A}$ is not a $\varphi$-badly approximable linear subspace. Then
there exist infinitely many $\bm{z}\in \Lambda \backslash \{0\}$  \ such that
$ {\rm dist} \,(
\bm{z} , \mathfrak{C})
 <\varphi (|\bm{z}|)
$
and then  
$\mathfrak{C}$ belongs to the set
$$
\{\mathfrak{C} =
\mathfrak{C}(\Theta) : 
\,\,\,   \Theta \,\,\text{satisfies (\ref{tt})},   \,\,\, 
$$
$$
\forall M \in \mathbb{N} \quad \exists \bm{z}\in \Lambda \backslash \{0\} : |\bm{z}| > M \quad \mbox{such that} \quad U_{\varphi(|\bm{z}|)} (\bm{z}) \cap \mathfrak{C} \neq \emptyset \} \subset
$$
 
 {
 
$$
\subset \bigcap_{M \in \mathbb{N}} \bigcup_{\substack{\bm{z} \in \Lambda\backslash \{0\} \\|\bm{z}| > M  \\ U_{\varphi(|\bm{z}|)} (\bm{z}) \cap \mathfrak{A} \neq \emptyset }} \widehat{L}_{\bm{z}^*} (R, {\varphi(|\bm{z}|)} )
= \bigcap_{M \in \mathbb{N}} \bigcup_{T\ge M}\,{G_T},
$$
where
$$
{G_T} = \bigcup_{\substack{\bm{z} \in \Lambda\backslash \{0\} \\|\bm{z}| =T  \\ U_{\varphi(T)} (\bm{z}) \cap \mathfrak{A} \neq \emptyset }} \widehat{L}_{\bm{z}^*} (R, {\varphi(T)} )
$$
In order to prove Theorem 2 we must show that
 $$ \mu\left(\bigcap_{M \in \mathbb{N}} \bigcup_{T\ge M}\,{G_T}\right) 
 =0, $$
and it is enough to show that the series 
\begin{equation}\label{conva}
\sum_{T=1}^\infty \mu\left({G_T}\right)
\end{equation}
converges.

By \eqref{varphi} for all $ T\in \mathbb{N}$ for $(a-c)\times c$-dimensional Lebesgue measure we have
$$
\mu 
\left({G_T}\right) \le
 R^{a+2}
   \dfrac{\varphi^{a-c}(T)}{T^{a-c}} K_T, 
$$
where $K_T = \# \{ \bm{z}  \in \Lambda: |\bm{z}| = T, \,\, \mathfrak{A} \cap U_{\varphi (|\bm{z}|)} (\bm{z}) \neq \emptyset \}$
and the convergence of the series (\ref{conva}) will follow from the convergence of the series
\begin{equation}\label{conva1}
\sum_{T=1}^\infty   \dfrac{\varphi^{a-c}(T)}{T^{a-c}} K_T.
\end{equation}
We should note that 
$$
\sum_{T/2< j\le T} K_j \le 
H_T = \# \{ \bm{z}  \in \Lambda : T/2<|\bm{z}| \le T, \mathfrak{A}\cap U_{\varphi (|\bm{z}|)} (\bm{z}) \neq \emptyset \}.
$$
We see that 
$$
\{ \bm{z}  \in \Lambda :  T/2< |\bm{z}| \le T, \,\, \mathfrak{A} \cap U_{\varphi (|\bm{z}|)} (\bm{z}) \neq \emptyset \}\subset
$$
$$
\subset
\{ \bm{z}  \in \Lambda :  |\bm{z}| \le T, \,\, \mathfrak{A}  \cap U_{\varphi (T/2)} (\bm{z}) \neq \emptyset \}\subset
\{ \bm{z}  \in \Lambda : \,\,\pmb{z}  \in 2\Pi_{T/2}   \},
$$
by monotonicity of $\varphi (\cdot)$.
 So
 $H_T \le N_{T/2}'$ by Lemma 2. But $N_{T/2}'  \ll \mu_T$
 by Remark to Lemma 1. Now 	we have
$$
 \sum_{j \le T} K_j \ll \mu_T+ \mu_{T/2}+ \mu_{T/2^2}+ \mu_{T/2^3}+ ... = M_T
.$$
We consider the partial sum of the series (\ref{conva1}) and apply partial summation to get
$$
\sum_{T=1}^{W+1}   \dfrac{\varphi^{a-c}(T)}{T^{a-c}} K_T=
$$
$$
=
 \sum_{T= 1}^{W+1} 
 \left(
 \dfrac{\varphi^{a-c}(T)}{T^{a-c}} - \dfrac{\varphi^{a-c}(T+1)}{(T+1)^{a-c}}
 \right)
  \cdot \sum_{j = 1}^T K_j + \dfrac{\varphi^{a-c}(W+1)}{(W+1)^{a-c}} \sum_{j= 1}^{W+1} K_j \ll
$$
 $$
 \ll
 \sum_{T= 1}^{W+1} 
 \left(
 \dfrac{\varphi^{a-c}(T)}{T^{a-c}} - \dfrac{\varphi^{a-c}(T+1)}{(T+1)^{a-c}}
 \right)
  \cdot M_T + \dfrac{\varphi^{a-c}(W+1)}{(W+1)^{a-c}}  M_{W+1}=
 $$
 $$
 =\sum_{T=1}^{W+1}   \dfrac{\varphi^{a-c}(T)}{T^{a-c}} \lambda_T.
 $$
 Here we use that $\frac{\varphi(T)}{T}$ decreases  and the equality (\ref{maa}).
By the condition of the Theorem \ref{th1} the series 
(\ref{series}) converges, and so the series  (\ref{conva1})  and (\ref{conva}) converge also. Theorem 2 is proved.$\Box$
}

\end{document}